\mag=\magstep1
\documentclass[reqno]{amsart}
\usepackage{graphicx}
\usepackage{amssymb}

\def\theorem#1{\bigbreak\noindent{\bf Theorem {#1}. }}

\def\lemma#1{\bigbreak\noindent{\bf Lemma {#1}. }}

\def\bgns{\it}

\def\fin{\rm\bigbreak}

\def\proof{\noindent{\it Proof. }}

\def\trans{\raise.8em\hbox{\sevenit t}}

\newcommand\al{\alpha}
\newcommand\be{\beta}

\newcommand\tht{\theta}

\newcommand\la{\lambda}

\newcommand\varf{\varphi}

\newcommand\w{\omega}

\newcommand\bZ{\mathbb{Z}}
\newcommand\tc{\tilde c}
\newcommand\td{\tilde d}

\newcommand\eqqx[1]{\begin{equation}#1\end{equation}}            
\newcommand\eqqs[1]{\begin{eqnarray}#1\end{eqnarray}}            

\newcommand\arrx[2]{\begin{array}{#1}#2\end{array}}

\newcommand\circled[1]{{\hskip3pt#1\hskip-.75em\raise.15ex\hbox{$\bigcirc$}}}

\newcommand\paren[1]{\left(#1\right)}
\newcommand\parenb[1]{\left[#1\right]}

\newcommand\dspl{\displaystyle}

\newcommand\vs{\vspace{1ex}}

\newcommand\arccosh{\operatorname{arccosh}}        
\newcommand\arcsinh{\operatorname{arcsinh}}        

\begin{document}

\title[Bipartite Chebyshev polynomials]
{\LARGE Bipartite Chebyshev polynomials\\ and elliptic integrals expressible\\ 
by elementary functions}
\author{KAZUTO ASAI}
\dedicatory{Center for Mathematical Sciences, University of Aizu, \vs\\
Aizu-Wakamatsu, Fukushima 965-8580, Japan \vs\\
{\rm e-mail: k-asai@u-aizu.ac.jp} \vs\\
Tel. 0242-37-2644 (Office), 0242-37-2752 (Fax)}
\address{Center for Mathematical Sciences, University of Aizu, Aizu-Wakamatsu, Fukushima 965-8580, Japan}
\email{k-asai@u-aizu.ac.jp}

\maketitle

\hrule \vs\vs
Abstract \\
The article is concerned with polynomials $g(x)$ whose graphs are ``partially packed'' between two horizontal tangent lines. We assume that most of the local maximum points of $g(x)$ are on the first horizontal line, and most of the local minimum points on the second horizontal line, except several ``exceptional'' maximum or minimum points, that locate above or under two lines, respectively. In addition, the degree of $g(x)$ is exactly the number of all extremum points $+1$. Then we call $g(x)$ a multipartite Chebyshev polynomial associated with the two lines. 

Under a certain condition, we show that $g(x)$ is expressed as a composition of the Chebyshev polynomial and a polynomial defined by the $x$-component data of the exceptional extremum points of $g(x)$ and the intersection points of $g(x)$ and the two lines. Especially, we study in detail bipartite Chebyshev polynomials, which has only one exceptional point, and treat a connection between such polynomials and elliptic integrals. 
\vskip2ex 
\hrule 
\bigskip

\section{Introduction}

The article is concerned with construction of real polynomials $g(x)$ whose graphs are partially packed between two horizontal tangent lines. Let $l_1,l_2$ be horizontal lines arranged downwards. We consider the case that most of the local maximum points of $g(x)$ are on $l_1$ and most of the local minimum points on $l_2$, except that several local maximum or minimum points, which we call the exceptional maximum or minimum points, are located above or under the two lines, respectively. Suppose the degree of $g(x)$ is exactly the number of all extremum points $+1$. Then we call $g(x)$ a multipartite Chebyshev polynomial associated with the two lines. The simplest case is that $g(x)$ has no exceptional extremum points, then $g(x)$ is essentially the Chebyshev polynomial of the first kind $T_n(x)$. For basics and several recent results of this polynomial, see \cite{cheby:theorie54, yusuf:cheby-poly17, mason:cheby02, erem-lem:extremal94, nash:cheby86, rivlin:cheby90, wat-zeit:min-poly93, zwillin:crcs11}. 

First of all, we consider the above-mentioned ``Chebyshev polynomial'' case. Let $m$ be a positive real constant. Let the horizontal tangent lines $l_1,l_2$ of $g(x)$ be $y=\pm m$ without loss of generality. Suppose $g(x)$ has positive degree $n$ and no exceptional extremum points. $g(x)-m$ has zeros at every local maximum point of $g(x)$, and also $g(x)+m$ has zeros at every local minimum point of $g(x)$. Additional two zeros of $g(x)-m$ or $g(x)+m$ exist, and we can set them $x=\pm a$ without loss of generality. Hence we have
\eqqx{n^2(g^2-m^2)=(x^2-a^2)g'^2\label{cheby-deq},}
and therefore, for $-a\le x\le a$, 
\eqqx{\int\frac{dg}{\sqrt{m^2-g^2}}=\pm n\int\frac{dx}{\sqrt{a^2-x^2}}.}
Letting $g=m\cos\varf$, $x=a\cos\tht$, and noting that $|x|\rightarrow a$ implies $|g|\rightarrow m$, we have
\eqqx{\arrx{c}{-\varf=\mp n\tht+k\pi.\qquad(k\in\bZ) \vs\\
\therefore\ \ g=\pm m\cos\paren{n\arccos\frac{x}{a}}=\pm mT_n\paren{\frac{x}{a}}.}}

We can treat similarly a multipartite Chebyshev polynomial $g(x)$ of positive degree $n$. Let the horizontal tangent lines $l_1,l_2$ of $g(x)$ be $y=\pm m$. Let $\al_1,\dots,\al_r$ be the $x$-components of the intersection (not tangent) points of the graph of $g(x)$ and $l_1$ or $l_2$, and let $\be_1,\dots,\be_{\ell}$ be the $x$-components of the exceptional extremum points of $g(x)$ outside the two lines. By definition, $r$ is even and $r=2\ell+2$. For convenience, we call the data $(\al_1,\dots,\al_r;\,\be_1,\dots,\be_{\ell})$ the outside data of $g(x)$, and set $p(x)=(x-\al_1)\dots(x-\al_r)$, $q(x)=(x-\be_1)\dots(x-\be_{\ell})$. We have
\eqqx{n^2(q(x))^2(g^2-m^2)=p(x)g'^2.\label{ode-for-g}}
Hence as above, for $x$ such that $p(x)<0$, 
\eqqx{\int\frac{dg}{\sqrt{m^2-g^2}}=\pm n\int\frac{q(x)}{\sqrt{-p(x)}}\,dx.
\label{int-for-g}}
In general, the RHS is a hyper-elliptic integral. Also, there is no assurance that \eqref{ode-for-g} or \eqref{int-for-g} has a polynomial solution $g(x)$. Now we assume that for some positive divisor $s$ of $n$, there exists a multipartite Chebyshev polynomial $u(x)$ of degree $s$ that shares the outside data with $g(x)$. Then $u$ satisfies the similar equation as \eqref{ode-for-g}:
\eqqx{s^2(q(x))^2(u^2-m^2)=p(x)u'^2.\label{ode-for-u}}
By \eqref{ode-for-g},\eqref{ode-for-u}, noting again that $x\rightarrow\al_i$ implies $|g|,|u|\rightarrow m$, and letting $g=m\cos\varf$, $u=m\cos\tht$, 
\eqqs{-\varf=&\dspl\int\frac{dg}{\sqrt{m^2-g^2}}=
\pm\frac{n}{s}\int\frac{du}{\sqrt{m^2-u^2}} \nonumber\vs\\
=&\dspl\mp\frac{n}{s}\tht+k\pi.\qquad(k\in\bZ) \nonumber\vs\\
\therefore\ \ g=&\pm m\cos\paren{\frac{n}{s}\arccos\frac{u}{m}}=
\pm mT_{\frac{n}{s}}\paren{\textstyle\frac{u}{m}}.
\label{super-ch-u}}

Equations \eqref{ode-for-g},\eqref{ode-for-u} are necessary conditions for $g$ and $u$, respectively, and they assure that the non-exceptional extremum points of them are located on $l_1,l_2$, but do not assure the total number of such points. In this case, however, as $u$ in \eqref{super-ch-u} is defined to be a multipartite Chebyshev polynomial, one can confirm that the total number of extremum points of $g$ is $(n/s)s-1=n-1$, which ensures $g$ to be a multipartite Chebyshev polynomial as desired. 

\theorem1\bgns Let $n$ be a positive integer and $s$ be a positive divisor of $n$. Let $g$ and $u$ be multipartite Chebyshev polynomials of degree $n$ and $s$, respectively, associated with the lines $y=\pm m$, and sharing a common outside data. Then we have $g=\pm mT_{\frac{n}{s}}\paren{\textstyle\frac{u}{m}}$. \fin

\bigskip

\section{Bipartite Chebyshev polynomials}

A multipartite Chebyshev polynomial $g(x)$ with a unique exceptional extremum point is called a bipartite Chebyshev polynomial. In this case there are four intersection points of $g(x)$ and two tangent lines $y=\pm m$. Suppose $g$ has positive degree $n$ and let the outside data of $g$ be $(\al_1,\al_2,\al_3,\al_4;0)$ $(\al_1>\al_2>0>\al_3>\al_4)$ without loss of generality. If the data satisfies symmetric condition: $\al_1=-\al_4$, $\al_2=-\al_3$, and $n$ is even, then $g$ is a special case of \eqref{super-ch-u} with quadratic $u$, say, $g=\pm mT_{\frac{n}{2}}\paren{\frac{2x^2-\al_1^2-\al_2^2}{\al_1^2-\al_2^2}}$. 

We now proceed to the case without symmetric condition. In accordance with the outside data, we have
\eqqx{n^2x^2(g^2-m^2)=p(x)g'^2.\label{g2-m2-p}}
Let $s$ be a positive divisor $s$ of $n$ (possibly $s=n$) such that there exists a bipartite Chebyshev polynomial $u$ of degree $s$ sharing a common outside data with $g$, then
\eqqx{s^2x^2(u^2-m^2)=p(x)u'^2.\label{u2-m2-p}}
By Theorem~1, $g$ is represented as $g=\pm mT_{\frac{n}{s}}\paren{\textstyle\frac{u}{m}}$. Thus we study the solution $u$ of equation \eqref{u2-m2-p}. Since $u'(0)=0$, we can set
\eqqx{u=a_0+a_2x^2+a_3x^3+\cdots+a_sx^s.\label{def-of-u}}
Dividing both sides of \eqref{u2-m2-p} by $x^2$ and differentiating with respect to $x$, we have
\eqqx{2s^2u=2u''\tilde p+u'{\tilde p}',\label{bi-ch-ode}}
where $\tilde p=p(x)/x^2$. It follows from \eqref{bi-ch-ode}, by setting $\tilde p=x^2+c_1x+c_2+c_3x^{-1}+c_4x^{-2}$ and comparing the coefficients of $x^k$, that 
\eqqx{a_k=\frac{1}{2(s^2-k^2)}\sum_{i=1}^4(k+i)(2k+i)c_ia_{k+i}\qquad(k=0,1,\dots,s-1),
\label{recur-bi-ch}}
where we promise $a_1=a_{s+1}=a_{s+2}=a_{s+3}=0$. By \eqref{recur-bi-ch}, the coefficients are determined step by step from $a_s$ to $a_3,a_2,a_0$ in the form:
\eqqx{a_k=F_k(c_1,c_2,c_3,c_4)a_s,\label{sol-rec-bi}}
where $F_k$ is a polynomial in $c_1,c_2,c_3,c_4$ with rational coefficients for every $k=0,2,3,$ $\dots,s$. The rest of the condition \eqref{u2-m2-p} is the equation \eqref{sol-rec-bi} for $k=1$:
\eqqx{F_1(c_1,c_2,c_3,c_4)=0,\label{condition-f1}}
and the $x^2$ terms of \eqref{u2-m2-p}:
\eqqx{s^2(a_0^2-m^2)=4a_2^2c_4\iff (s^2F_0^2-4c_4F_2^2)a_s^2=s^2m^2
.\label{m-equation}}

\lemma1\bgns For every $k=0,1,\dots,s$, the polynomial $F_k=F_k(c_1,c_2,c_3,c_4)$ has the following properties:
\begin{itemize}
\item[(i)] The coefficient of every term of $F_k$ is positive.
\item[(ii)] $F_{s-k}$ consists of the terms of type $c_{i_1}c_{i_2}\dots c_{i_r}$, such that $i_1+i_2+\cdots+i_r=k$.
\item[(iii)] Every possible term of type $c_{i_1}c_{i_2}\dots c_{i_r}$ with $i_1+i_2+\cdots+i_r=k$ appears in $F_{s-k}$ except that $F_0$ does not contain the term $c_1^s$. 
\end{itemize}\fin

\proof For convenience, write $c_{i_1}c_{i_2}\dots c_{i_r}=c_{\la}$, where $\la=(\la_1,\dots,\la_r)$ is the weakly decreasing rearrangement of $i_1,\dots,i_r$. Then $\la$ is a partition of $k$ (denoted by $\la\vdash k$), $k=i_1+\cdots+i_r$. Each $\la_i$ is called a part of $\la$, and $r$ is called the length of $\la$ denoted by $\ell(\la)$. There is a unique partition of zero, denoted by $\varnothing$, consisting of no parts, having length 0, and we promise that $c_{\varnothing}=1$. There are no partitions of negative integers. 

By the fact that all coefficients of the linear recurrence relation \eqref{recur-bi-ch} are positive and by \eqref{sol-rec-bi}, we see (i). Next, we prove (ii) and (iii) by induction on $k$. When $k=0$, we have $F_s=1$ which satisfies the proposition. For the induction step, suppose, for every nonnegative integer $j<k$, 
\eqqx{F_{s-j}=\sum_{\la\vdash j;\,\la_1\le4}m_{\la}c_{\la}\label{f-exprsdbymc}}
with some positive coefficients $m_{\la}$. Then by \eqref{recur-bi-ch}, \eqref{sol-rec-bi}, putting $\frac{(s-k+i)(2s-2k+i)}{2k(2s-k)}=m_i^{(k)}$, 
\eqqx{\arrx{rl}{F_{s-k}&=\dspl\sum_{i=1}^4m_i^{(k)}c_iF_{s-k+i}=\sum_{i=1}^4
\sum_{\la\vdash k-i;\,\la_1\le4}m_i^{(k)}m_{\la}c_ic_{\la} \\&=
\dspl\sum_{\la\vdash k;\,\la_1\le4}m_{\la}c_{\la}.}\label{recur-fs}}
Here, as every partition $\la$ of $k$ with $\la_1\le4$ is composed of the parts $\la_j$ such that $1\le\la_j\le4$, we can reduce the partition $\la$ to a partition of $k-i$ by subtracting some part $1\le i\le4$. Therefore, together with the positivity of the coefficients $m_i^{(k)}m_{\la}$, every coefficient $m_{\la}$ is shown to be positive. Only the case $F_0$, however, is represented as
\eqqx{F_0=\dspl\sum_{i=2}^4m_i^{(s)}c_iF_i=\sum_{i=2}^4
\sum_{\la\vdash s-i;\,\la_1\le4}m_i^{(s)}m_{\la}c_ic_{\la},} 
because $a_1$ is defined to be 0. From this, it follows that $F_0$ does not contain only the term $c_{(1,\dots,1)}=c_1^s$. \qed\fin

Reviewing \eqref{recur-fs}, we have $m_{\la}=\sum_{i=1}^4m_i^{(k)}m_{\la\ominus(i)}$ for $\la\vdash k$, where $\la\ominus(i)$ denotes a partition obtained by removing a part of size $i$ from $\la$, and set $m_{\la\ominus(i)}=0$ in case that $\la$ has no parts of size $i$. Iterating this induction, we can represent the coefficient $m_{\la}$ as follows. 

\lemma2\bgns Let $S(\la)$ denote the set of all distinct permutations $i=(i_1,\dots,i_r)$ of $\la$. If $\la$ is a partition of the integer less than $s$, then 
\eqqx{m_{\la}=\sum_{i\in S(\la)}\prod_{t=1}^rm_{i_t}^{(i_1+\cdots+i_t)},\label{m-la-formula}}
while if $\la\vdash s$, $S(\la)$ in this formula should be replaced with $S'(\la)$, the set of all distinct permutations $(i_1,\dots,i_r)$ of $\la$ where $i_r>1$. \fin

On the basis of the above-mentioned method, we construct the polynomial $u$ as described in the procedure below:
\begin{itemize}
\item[(i)] Choose an arbitrary positive constant $m$, and arbitrary real coefficients $c_2,c_3,c_4$ of $p$ such that $c_2<0<c_4$. 
\item[(ii)] Represent all coefficients $a_k$ of $u$ as a polynomial in $c_1,\dots,c_4$ and $a_s$ by using \eqref{sol-rec-bi}, \eqref{f-exprsdbymc} and \eqref{m-la-formula}. 
\item[(iii)] By \eqref{condition-f1}, determine $c_1$ depending on $c_2,c_3,c_4$. 
\item[(iv)] Determine $a_k$ for all $k=0,2,3,\dots,s-1$. 
\item[(v)] Also determine $a_s$ by \eqref{m-equation}. 
\end{itemize}
Unfortunately, in the steps (iii) and (v), $c_1$ is possibly an imaginary number, and $a_s$ may be also imaginary, or even does not exist if $s^2F_0^2-4c_4F_2^2=0$. Now we avoid those difficulties by handling the coefficients $c_3,c_4$ ($c_2<0$ is fixed). Letting $c_3,c_4=0$, equation \eqref{u2-m2-p} becomes 
\eqqx{s^2x^2(u^2-m^2)=x^2(x^2+c_1x+c_2)u'^2,\label{limit-bipar}}
which is clearly reduced to \eqref{cheby-deq}, and we have $u=\pm mT_s\paren{\frac{x-b}{a}}$. But by the continuity of all equations in the above process (i)---(v), we can show that for $c_3$ and $c_4$ sufficiently close to 0, $c_1$ and $a_s$ are determined to be real numbers as follows. Equation \eqref{limit-bipar} has $s-1$ possible values of $c_1$ for given $c_2$, because the condition \eqref{condition-f1}, i.e., $u'(0)=0$ obliges one of the extremum points of $u$ to be located on the $y$-axis, which allows $s-1$ possible graphs of $u=u_1,\dots,u_{s-1}$. While by Lemma~1, $F_1$ has a leading term $c_1^{s-1}$ with respect to $c_1$, and this has $s-1$ distinct real roots $c_1$ for $c_3=c_4=0$. Hence, again for $c_3$ and $c_4$ $(c_4>0)$ sufficiently close to 0, the outside data of $u$ satisfies $\al_1>\al_2>0>\al_3>\al_4$, and by the continuity of $F_1$, we have $s-1$ distinct real values of $c_1$, and therefore $s-1$ different $u$'s (also possible two $\pm a_s$'s), that are obtained graphically from $u_1,\dots,u_{s-1}$ moving slightly the maximum or minimum point on $y$-axis upward or downward, respectively. Reducing the expression $g=\pm mT_{\frac{n}{s}}\paren{\frac{u}{m}}$ to $\pm mT_{\frac{n}{s}}(u)$, we have the following. 

\theorem2\bgns For every positive integer $n$ and every positive divisor $s$ of $n$, there exists a bipartite Chebyshev polynomial of degree $n$ of the form $g=\pm mT_{\frac{n}{s}}(u)$ associated with the lines $y=\pm m$, where $u$ is also a bipartite Chebyshev polynomial of degree $s$ associated with the lines $y=\pm1$ obtained by moving an arbitrary one of the maximum or minimum points of $\pm T_s\paren{\frac{x-b}{a}}$ upward or downward, respectively. \fin

\bigskip

\section{Elliptic integrals}

It is known that, in general, elliptic integrals are not representable in terms of elementary functions. In this section, however, for a special kind of elliptic integral, we show that there are infinitely many representable cases, which are represented as compositions of the inverse trigonometric/hyperbolic functions and polynomials or compositions of the logarithm and polynomials. For general properties of elliptic integrals and elliptic functions, see e.g. \cite{hanc:ellip17, pra-sol:ellip97}.

Focusing again on the integral solution of \eqref{u2-m2-p}, a special case of the solution \eqref{int-for-g} (with $g$, $n$ replaced by $u$, $s$, respectively), for an arbitrary monic real quartic polynomial $p(x)$, 
\eqqx{\arrx{rl}{\dspl\int\frac{x}{\sqrt{-p(x)}}\,dx=\pm\frac{1}{s}\arccos\frac{u}{m}+C
&\paren{\arrx{l}{\mbox{in each interval such that $p(x)<0$} \\ \mbox{and $-m<u(x)<m$}}} 
\vs\\
\dspl\int\frac{x}{\sqrt{p(x)}}\,dx=\pm\frac{1}{s}\arccosh\frac{|u|}{m}+C
&(\mbox{in each interval such that $p(x)>0$}).}\label{bipartite-int}}
Hence, if $u$ is a polynomial (which should be of degree $s$), then the elliptic integrals on the LHS are represented as elementary functions. According to Section~2, for every positive integer $s$, there exists a polynomial solution $u$ of \eqref{u2-m2-p} of degree $s$, under suitable condition for $p(x)$, and in this case we obtain a bipartite Chebyshev polynomial $u$, which satisfies both of \eqref{bipartite-int}. For more general $p(x)$, if we ignore the outside data, and only suppose $x^2\nmid p(x)$, then the polynomial solution $u$ of \eqref{u2-m2-p} satisfies $u'(0)=0$, and therefore \eqref{condition-f1} and \eqref{m-equation} are satisfied, and conversely, they assure the existence of the polynomial solution. 

For the first formula of \eqref{bipartite-int}, the domain $D$ defined by $p(x)<0$ is divided by the condition $-m<u(x)<m$ into small intervals to avoid the singularities of $\arccos x$, $\pm1$. If we consider the RHS in $D$, ignoring the integral constant and $-$ sign at the beginning, the graph consists of zigzag-like lines packed between the horizontal lines $y=\frac{\pi}{s}$ and $x$-axis, and a small piece of the graph between cusps matches the LHS with the appropriate sign. To complete the formula in $D$, we should ``join'' the pieces of the RHS changing the sings and shifting the integral constant by an integral multiple of $\frac{\pi}{s}$. 

For the second formula of \eqref{bipartite-int}, we see that $u(x)>m$ or $u(x)<-m$ is not necessary to add to the condition $p(x)>0$. 

Suppose that $u$ is a bipartite Chebyshev polynomial of degree $s$ with the outside data $(\al_1,\al_2,\al_3,\al_4;0)$ $(\al_1>\al_2>0>\al_3>\al_4)$, then $D$ is composed of the intervals $(\al_4,\al_3)$, $(\al_2,\al_1)$. Let the unique exceptional extremum point of $u(x)$ be the $k$-th extremum point counted from the right. From the above consideration of the zigzag-like graph of the first formula, it follows that
\eqqx{\w_1=\int_{\al_2}^{\al_1}\frac{x}{\sqrt{-p(x)}}\,dx=\frac{k}{s}\pi;\quad
\w_2=\int_{\al_4}^{\al_3}\frac{x}{\sqrt{-p(x)}}\,dx=\frac{k-s}{s}\pi;\quad\w_1-\w_2=\pi.}

Forgetting the connection with multipartite Chebyshev polynomials, we can also deal with the equation similar to \eqref{g2-m2-p} for a polynomial $g$ of positive degree $n$:
\eqqx{n^2x^2(g^2+m^2)=p(x)g'^2,\label{g2+m2-p}}
where $p(x)$ is a monic real quartic polynomial satisfying $x^2\nmid p(x)$. In the same way as in Section~2, for some positive divisor $s$ of $n$, suppose there exists a polynomial $u$ of degree $s$ such that
\eqqx{s^2x^2(u^2+m^2)=p(x)u'^2.\label{u2+m2-p}}
Letting $g=m\sinh\varf$, $u=m\sinh\tht$ temporarily in the complex domain (to determine the constant of integration), and noting that $p(x)\rightarrow0$ implies $g,u\rightarrow\pm mi$, we have 
\eqqs{\varf=\dspl\int\frac{dg}{\sqrt{m^2+g^2}}=
\pm\frac{n}{s}\int\frac{du}{\sqrt{m^2+u^2}}=\dspl\pm\frac{n}{s}\tht.\nonumber\vs\\
\therefore\ \ g=\pm m\sinh\paren{\frac{n}{s}\arcsinh\frac{u}{m}}\equiv
\pm m\tilde T_{\frac{n}{s}}\paren{\textstyle\frac{u}{m}}.
\label{super-shash-u}}
Here, if $g,u\rightarrow\pm mi$, then $\varf\rightarrow(\frac{\pi}{2}+k\pi)i$, $\tht\rightarrow(\frac{\pi}{2}+k'\pi)i$ $(k,k'\in\bZ)$, but this implies that $n/s$ is odd, which we now assume, and the (real) constant of integration vanishes. One sees that $\tilde T_n(x)=\sinh\paren{n\arcsinh x}$ is monotonously increasing, and a polynomial of degree $n$ iff $n$ is an odd positive integer. 

Returning to \eqref{u2+m2-p}, a change from \eqref{u2-m2-p} to \eqref{u2+m2-p} does not have an effect on the condition \eqref{condition-f1}, while it turns \eqref{m-equation} into
\eqqx{s^2(a_0^2+m^2)=4a_2^2c_4\iff (4c_4F_2^2-s^2F_0^2)a_s^2=s^2m^2.
\label{-m-equation}}
For convenience, set $d=s^2F_0^2-4c_4F_2^2$. We see the signature of $d$ discriminates the existence of $a_s$ in \eqref{m-equation} or \eqref{-m-equation}. On the other hand, the solution of \eqref{g2+m2-p} is represented as
\eqqx{\int\frac{x}{\sqrt{p(x)}}\,dx=\pm\frac{1}{s}\arcsinh\frac{u}{m}+C.
\label{hyper-bip-int}}

For the special case $m=0$, one sees that 
\eqqx{\int\frac{x}{\sqrt{p(x)}}\,dx=\pm\frac{1}{s}\log|u|+C.
\label{log-int}}

\theorem3\bgns Let $s$ be a positive integer. For a real quartic polynomial $p(x)=x^4+c_1x^3+c_2x^2+c_3x+c_4$ such that $x^2\nmid p(x)$, the elliptic integral 
\eqqx{\int\frac{x}{\sqrt{\pm p(x)}}\,dx\label{ellip-int}}
is represented as either \eqref{bipartite-int} or \eqref{hyper-bip-int} or \eqref{log-int} subject to $d>0$ or $d<0$ or $d=0$, respectively, where $u$ is a polynomial of degree $s$ determined by \eqref{def-of-u}, \eqref{sol-rec-bi}, \eqref{f-exprsdbymc}, \eqref{m-la-formula} and \eqref{m-equation} or \eqref{-m-equation}, if and only if \eqref{condition-f1} is satisfied. 
\fin

Lastly, we study the set of quartic polynomials $p(x)$ such that \eqref{ellip-int} is representable as in Theorem~3, depending on the degree $s$ of $u$. For even $s$, the polynomial $F_1(c_1,c_2,c_3,$ $c_4)$ has the leading term $m_{(1,\dots,1)}c_1^{s-1}$ with respect to $c_1$, and thus (as $s-1$ is odd), for arbitrary $c_2,c_3,c_4$, we can find the real solution $c_1$ to \eqref{condition-f1}. Therefore for given $c_2,c_3,c_4$, we have $p(x)$ suitable for our representation as in Theorem~3. For $s=4k-1$, $F_1$ has the leading term $m_{(2,\dots,2)}c_2^{2k-1}$ with respect to $c_2$, and therefore for arbitrary $c_1,c_3,c_4$, we have $c_2$ and $p(x)$ as desired. For $s=8k-3$, $F_1$ has the leading term $m_{(4,\dots,4)}c_4^{2k-1}$ with respect to $c_4$, and so for arbitrary $c_1,c_2,c_3$, we have $c_4$ and $p(x)$ similarly. The rest case $s=8k+1$ gives no clear assurance for $p(x)$. 

\theorem4\bgns If $s$ is not congruent to 1 modulo 8, then for some $i_1,i_2,i_3$ $(1\le i_1<i_2<i_3\le4)$, given $c_{i_1},c_{i_2},c_{i_3}$, we obtain $p(x)$ suitable for representation of \eqref{ellip-int} as in Theorem~3 by a polynomial $u$ of degree $s$. 
\fin

\bigskip

\bibliography{../rrefs/refchebyshev}
\bibliographystyle{plain}

\nocite{*}

\end{document}